\documentclass[11pt,a4paper,reqno]{amsart}
\usepackage[T1]{fontenc}
\usepackage{lmodern}
\usepackage{mathtools,amssymb}
\usepackage[margin=28mm]{geometry}
\usepackage{microtype,etoolbox}
\usepackage[hidelinks,bookmarksnumbered]{hyperref}
\usepackage{xcolor}
\usepackage{float}
\AtBeginDocument{%
  \begingroup\footnotesize\fontencoding{OT1}\selectfont\endgroup
}

\numberwithin{equation}{section}
\newtheorem{theorem}{Theorem}[section]
\newtheorem{proposition}[theorem]{Proposition}
\newtheorem{lemma}[theorem]{Lemma}
\newtheorem{corollary}[theorem]{Corollary}
\theoremstyle{definition}
\newtheorem{definition}[theorem]{Definition}

\theoremstyle{remark}
\newtheorem{remark}[theorem]{Remark}
\newtheorem{example}[theorem]{Example}

\BeforeBeginEnvironment{theorem}{\begin{samepage}}
\AfterEndEnvironment{theorem}{\end{samepage}\nopagebreak[3]}
\BeforeBeginEnvironment{proposition}{\begin{samepage}}
\AfterEndEnvironment{proposition}{\end{samepage}\nopagebreak[3]}
\BeforeBeginEnvironment{lemma}{\begin{samepage}}
\AfterEndEnvironment{lemma}{\end{samepage}\nopagebreak[3]}
\BeforeBeginEnvironment{corollary}{\begin{samepage}}
\AfterEndEnvironment{corollary}{\end{samepage}\nopagebreak[3]}
\newcommand{\R}{\mathbb R}

\newcommand{\C}{\mathbb C}
\newcommand{\Z}{\mathbb Z}
\newcommand{\spec}{\operatorname{spec}}
\newcommand{\Span}{\operatorname{span}}
\newcommand{\Tr}{\operatorname{Tr}}
\newcommand{\ip}[2]{\langle #1,#2\rangle}
\newcommand{\norm}[1]{\lVert #1\rVert}
\newcommand{\diag}{\operatorname{diag}}

\allowdisplaybreaks[1]
\makeatletter
\patchcmd{\ps@headings}
  {\normalfont\scriptsize}
  {\normalfont\fontsize{7}{9}\selectfont}{}{}
\patchcmd{\ps@headings}
  {\normalfont\scriptsize}
  {\normalfont\fontsize{7}{9}\selectfont}{}{}
\makeatother
\title{Comparing Lower and Upper Spectral Gaps of Nonnegative Symmetric Matrices
}
\author{Yanlong Ding}
\address{\textnormal{Yanlong Ding:}\endgraf\nobreak\indent School of Mathematical Sciences, University of Science and Technology of China, Hefei 230026, China.}
\email{dylustc@mail.ustc.edu.cn}

\author{Chunyang Hu}
\address{\textnormal{Chunyang Hu:}\endgraf\nobreak\indent Yau Mathematical Sciences Center, Tsinghua University, Beijing, 100084, China}
\email{chunyanghu@mail.tsinghua.edu.cn}

\author{Shiping Liu}
\address{\textnormal{Shiping Liu:}\endgraf\nobreak\indent School of Mathematical Sciences, University of Science and Technology of China, Hefei 230026, China.}
\email{spliu@ustc.edu.cn}

\author{Chiyu Zhou}
\address{\textnormal{Chiyu Zhou:}\endgraf\nobreak\indent School of Mathematical Sciences, University of Science and Technology of China, Hefei 230026, China.}
\email{dovong@mail.ustc.edu.cn}

\date{7 November 2026}
\begin{document}

\begin{abstract}
We establish sharp estimates comparing the lower and upper spectral gaps of nonnegative symmetric matrices.
Combining these comparisons with a multiplicity argument in the spirit of Li and Liu, we improve and extend earlier spectral gap comparison bounds of Saha and Li and Liu, which provide quantitative refinements of an observation of Breuillard, Green, Guralnick and Tao relating combinatorial expansion to two-sided spectral expansion. Our graph estimates hold for nonbipartite even-walk-regular graphs, which covers the classes of graphs, including vertex-transitive graphs, (twisted) Cayley graphs, and (twisted) Cayley sum graphs, investigated in previous work.

Under an additional eigenmap assumption, we also establish quantitative two-sided comparisons between the lower and upper spectral gaps. In particular, these comparisons apply to walk-regular graphs.
\end{abstract}
\maketitle
\markboth{\shorttitle}{\shorttitle}
\section{Introduction}

\subsection{Background and related works} Let $G=(V,E)$ be a finite connected simple graph on $n\ge2$ vertices,
with adjacency matrix $A$ and degree matrix $D$. Its normalized
adjacency matrix
\[
 P=D^{-1}A
\]
is similar to the symmetric nonnegative matrix $D^{-1/2}AD^{-1/2}$.
Write its eigenvalues, counted with multiplicity, as
\[
 1=\mu_1>\mu_2\ge\cdots\ge\mu_n\ge-1.
\]
The eigenvalue $1$ is simple because $G$ is connected, and
$\mu_n=-1$ if and only if $G$ is bipartite.
For $1\le k\le n$, we call $1+\mu_{n-k+1}$ the
\emph{$k$-th lower spectral gap} and $1-\mu_k$ the
\emph{$k$-th upper spectral gap}. In this article, we study comparisons between lower and upper spectral gaps.

Eigenvalue estimates on graphs have attracted a lot of attention. In particular, bounds in terms of the \emph{Cheeger constant} $h$ have been extensively studied \cite{AM85,BKW15,chung97,Dodziuk84,LS88,LOT,Miclo08,Miclo15}. Estimates involving the \emph{vertex isoperimetric constant} $h_{\mathrm{out}}$ were first investigated by Alon~\cite{Alon} in 1986, motivated by Cheeger's work on Riemannian manifolds~\cite{Cheeger}.

Breuillard, Green, Guralnick and Tao~\cite[Proposition~E.1]{BGGT}
showed that, for nonbipartite Cayley graphs of uniformly
bounded degree, a uniform positive lower bound on
$h_{\mathrm{out}}$ ensures that all normalized adjacency eigenvalues
other than the trivial eigenvalue $\mu_1=1$ are uniformly bounded
away from both $1$ and $-1$.

Quantitative versions of this observation were subsequently developed
by Biswas~\cite{Biswas}, Biswas and Saha~\cite{BSsum},
Moorman, Ralli and Tetali~\cite{MRT}, Saha \cite{Saha1} and Hu and Liu~\cite{HL},
yielding bounds on the nontrivial eigenvalues $\{\mu_2,\ldots,\mu_n\}$ in terms of $h_{\mathrm{out}}$ or $h$ and the degree $d$.
In particular, Hu and Liu~\cite{HL} obtained a lower bound of order
$h_{\mathrm{out}}^2/d$ for $1+\mu_n$ on nonbipartite
$d$-regular vertex-transitive graphs.
Li and Liu~\cite{LL} established quantitative comparisons between
$1+\mu_{n-1}$ and $1-\mu_2$ for connected graphs and derived
lower bounds for $1+\mu_n$ in terms of $1-\mu_2$ and vertex degree on nonbipartite vertex-transitive graphs.
Saha~\cite{Saha} proved a linear lower bound for $1+\mu_n$
in terms of $(1-\mu_2)/d$ on nonbipartite vertex-transitive graphs, (twisted) Cayley graphs and (twisted) Cayley sum graphs.

In our terminology, Saha's result~\cite{Saha} compares
the \emph{first} lower spectral gap \(1+\mu_n\) with the
\emph{second} upper spectral gap \(1-\mu_2\) for certain
classes of graphs. As pointed out by Li and
Liu~\cite[Remark~1]{LL}, the higher-order dual Cheeger
inequality~\cite{LIU2015} provides a lower bound for the
\emph{second} lower spectral gap \(1+\mu_{n-1}\) in terms
of the second upper spectral gap \(1-\mu_2\) for general graphs. More generally, for every \(1\leq k\leq n\), the higher-order
dual Cheeger inequality~\cite[Theorem~1.2]{LIU2015} gives
\[
    1+\mu_{n-k+1}
    \geq C(k)\bigl(1-\overline{h}(k)\bigr)^2,
\]
where \(\overline{h}(k)\) is the \(k\)-way dual Cheeger
constant defined in~\cite[(8)]{LIU2015}, and \(C(k)>0\)
depends only on \(k\). By~\cite[Proposition~3.1]{LIU2015},
\[
    1-\overline{h}(k)\geq h(k),
\]
where \(h(k)\) is the \(k\)-way Cheeger constant defined
in~\cite[(4)]{LIU2015}. The easy direction of the higher-order
Cheeger inequality~\cite{LOT}
(see also~\cite[Theorem~1.1]{LIU2015}) gives
\[
    1-\mu_k\leq 2h(k).
\]
Combining these inequalities, we obtain
\[
    1+\mu_{n-k+1}
    \geq C(k)\bigl(1-\overline{h}(k)\bigr)^2
    \geq C(k)h(k)^2
    \geq \frac{C(k)}{4}(1-\mu_k)^2.
\]
This motivates comparisons between the \(k\)-th lower
and upper spectral gaps.

In this paper, we improve this inequality and establish
sharp comparisons for both individual spectral gaps and
sums of spectral gaps over consecutive indices. Our results
hold for arbitrary nonnegative symmetric matrices.

To compare the first lower spectral gap \(1+\mu_n\) with
the second upper spectral gap \(1-\mu_2\), one can use an
argument based on the multiplicity of the smallest
eigenvalue, as in~\cite{LL}. If \(\mu_n=\mu_{n-1}\), any
lower bound for \(1+\mu_{n-1}\) immediately yields the same
bound for \(1+\mu_n\). For nonbipartite graphs, the remaining
case, in which \(\mu_n\) is simple, can be treated under
additional structural assumptions, such as vertex
transitivity. 

Applying our sharp comparisons between the \(k\)-th lower
and upper spectral gaps, 
we obtain improved lower bounds for \(1+\mu_n\) in terms
of \(1-\mu_2\) for nonbipartite even-walk-regular graphs.
The class of even-walk-regular graphs includes the graph
classes studied in previous work: vertex-transitive graphs,
(twisted) Cayley graphs,
and (twisted) Cayley sum graphs. Moreover, for walk-regular
graphs, we obtain quantitative comparisons in both
directions between the \(k\)-th lower and upper spectral gaps.

\subsection{Comparison of lower and upper spectral gaps}
\label{intro-sec:comparison}
Let $N\in\R^{n\times n}$ be a symmetric matrix with nonnegative
entries, that is, $N_{ij}\ge0$ for all $i,j$.
List its eigenvalues, with multiplicity, as
\begin{equation}\label{eq:lambda}
 \lambda_1\ge\lambda_2\ge\cdots\ge\lambda_n,
 \qquad \rho=\norm N=\lambda_1.
\end{equation}
Here $\norm \cdot$ denotes the operator norm and $\rho$ is the spectral radius.
The corresponding $k$-th lower and upper spectral gaps
are $\rho+\lambda_{n-k+1}$ and $\rho-\lambda_k$, respectively. The following theorem compares these gaps at each index.

\begin{theorem}[Lower and upper spectral gap comparisons]\label{ed-thm:main}
Let $N\in\R^{n\times n}$ be a symmetric matrix with nonnegative
entries, with eigenvalues $\lambda_1\ge\cdots\ge\lambda_n$
and spectral radius $\rho=\norm N=\lambda_1$.
For every $1\le m\le n$, we have
\begin{equation}\label{ed-eq:gap-scaled}
 \rho+\lambda_{n-m+1}
 \ge\frac{\rho-\lambda_m}{4}
\end{equation}
and
\begin{equation}\label{eq:prefix-gaps}
   \sum_{j=n-m+1}^n\left(\rho+\lambda_j\right)\ge\frac12\sum_{j=1}^m\left(\rho-\lambda_j\right).
\end{equation}
\end{theorem}
The theorem applies in particular to adjacency matrices of finite
simple graphs. Applying it to $D^{-1/2}AD^{-1/2}$ also gives
the corresponding inequalities for normalized adjacency
eigenvalues, with $\rho=1$.

Proposition~\ref{ed-prop:sharp-gap} shows that the factor $1/4$ in (\ref{ed-eq:gap-scaled}) is optimal for every fixed $m\ge2$, even among connected nonbipartite vertex-transitive simple graphs. And Example~\ref{ex:prefix-sharp} shows that the factor $1/2$ in the partial-sum inequality (\ref{eq:prefix-gaps}) is optimal for every
fixed $m\ge2$.

\subsection{Spectral symmetry and gap bounds for walk-regular graphs}
\label{intro-sec:eigenmap}
\label{intro-sec:ewr}
A graph $G=(V,E)$ with adjacency matrix $A$ is
\emph{even-walk-regular} if
\[
 (A^{2\ell})_{xx}=(A^{2\ell})_{yy}
 \qquad(x,y\in V,\ \ell\in\Z_{\ge0}).
\]
If, more strongly, $(A^\ell)_{xx}$ is independent of $x$ for every
$\ell\in\Z_{\ge0}$, then $G$ is \emph{walk-regular}.
Godsil and McKay~\cite{GM} introduced walk-regular graphs
in their study of cospectral graphs. Note that both conditions imply regularity, since $(A^2)_{xx}=\deg(x)$.

Our next result compares the first lower spectral gap with the
second upper spectral gap on even-walk-regular graphs.

\begin{theorem}\label{thm:ewr-gap}
For every finite connected nonbipartite simple even-walk-regular graph $G$ of degree $d$ with normalized adjacency eigenvalues $\{\mu_i\}$, we have
\begin{equation}\label{eq:ewr-gap}
 1+\mu_n\ge\min\left\{\frac2d,\frac{1-\mu_2}{4}\right\}
 \ge\frac{1-\mu_2}{2d}.
\end{equation}
\end{theorem}
Saha \cite[Theorem 1.1]{Saha} proved $1+\mu_n \ge (1-\mu_2)/(50000d)$ for connected nonbipartite Cayley, Cayley sum, twisted Cayley, twisted Cayley sum, and vertex-transitive graphs. Proposition~\ref{cg-prop:fiveEWR} in Appendix~\ref{app:five-classes} shows that every graph in these classes is even-walk-regular.
Thus Theorem~\ref{thm:ewr-gap} strengthens his bound on these simple graphs and extends it to all connected nonbipartite simple even-walk-regular graphs. For walk-regular graphs, we also obtain quantitative two-sided
comparisons between the $j$-th lower and upper spectral gaps. In addition, in Section \ref{ed-sec:sharp} we provide a sharp example: $G_m=C_m \lor C_m$, i.e. the join of two cycles $C_m$, to show that the estimate for the first lower spectral gap $1+\mu_n$ in terms of $h_{out}$ for a nonbipartite vertex-transitive graph, which is established by Hu and Liu \cite[Theorem 1.3]{HL}, is optimal in the order of $d$.

\begin{theorem}\label{rf-thm:rootgaps}
Let $G$ be a finite connected walk-regular simple graph with $n\ge2$ vertices, and $\{\mu_i\}$ be its normalized adjacency eigenvalues.
Then, for every $1\le j\le n$,
\begin{equation}\label{eq:wr-angle}
 |\arccos\mu_j+\arccos\mu_{n-j+1}-\pi|\le\arccos(-\mu_n),
\end{equation}
and
\begin{equation}\label{rf-eq:sqrtgaps}
 |\sqrt{1-\mu_j}-\sqrt{1+\mu_{n-j+1}}|
 \le\sqrt{1+\mu_n},\qquad
 |\mu_j+\mu_{n-j+1}|\le\sqrt{2(1+\mu_n)}.
\end{equation}

\end{theorem}
\begin{figure}[H]
    \centering
    \includegraphics[width=0.5\linewidth]{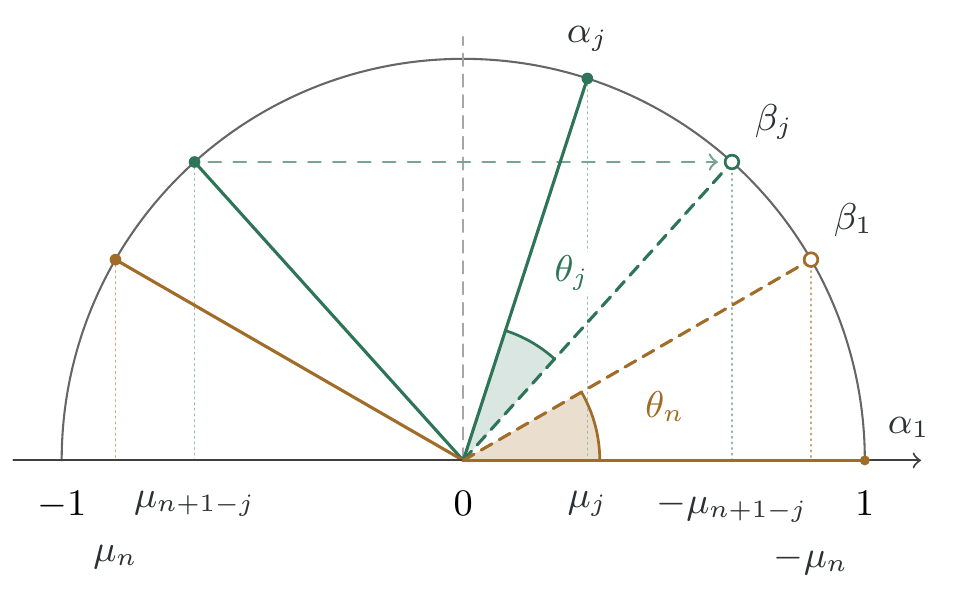}
  \caption{Angular reflection bound $\theta_j\le\theta_1=\theta_n$.
Solid rays make angles $\arccos\mu_j$ with the positive horizontal
axis; dashed rays are their reflections across the vertical axis.}
    \label{fig:placeholder}
\end{figure}
Using the identity $\arccos(-t)=\pi-\arccos t$, we may rewrite (1.5) as
\[
\left|\arccos\mu_j-\arccos(-\mu_{n+1-j})\right|
\le
\left|\arccos\mu_1-\arccos(-\mu_n)\right|,
\qquad 1\le j\le n.
\]
Define
\[
\theta_j:=
\left|\arccos\mu_j-\arccos(-\mu_{n+1-j})\right|,
\]
the angle between the rays through $\alpha_j$ and $\beta_j$.
Since $\mu_1=1$, we have
$\theta_1=\theta_n=\arccos(-\mu_n)$.
Thus, (1.5) states that every angle $\theta_j$ is bounded by the
endpoint angle:
\[
\theta_j\le\theta_1=\theta_n,
\qquad 1\le j\le n.
\]

\subsection{The structure of the paper} In Section \ref{section:preliminary}, we present basic definitions, notations and useful lemmas for the later discussions and proofs. In Section \ref{ed-sec:top}, we use the sphere case of Dold's theorem and Ky Fan's principle to compare lower and upper spectral gaps of symmetric nonnegative matrices, proving Theorem~\ref{ed-thm:main}. In Section \ref{rf-sec:angular}, we use eigenmaps to obtain two-sided spectral comparisons for symmetric nonnegative matrices and apply them to walk-regular and even-walk-regular graphs, proving Theorems~\ref{rf-thm:rootgaps} and~\ref{thm:ewr-gap}. In Section \ref{ed-sec:sharp}, we present various examples to show sharpness of the comparisons between eigenvalues. Appendix~\ref{app:five-classes} proves that the five graph classes considered by Saha are even-walk-regular.

\subsection{Disclosure of AI assistance}
Generative artificial-intelligence tools were used to develop the
proofs and assist with exposition. The authors take responsibility
for the mathematical content and its verification.

\section{Preliminary}\label{section:preliminary}
\subsection{Graph}
Let $G=(V,E)$ be a finite graph with vertex set $V$ and edge set $E\subseteq\{\{x,y\}:x\neq y\in V\}$. For two vertices $x,y\in V$, $x$ and $y$ are adjacent, if $\{x,y\}\in E$ and we denote as $x\sim y$. For a vertex $x\in V$, the degree of $x$ is $\vert \{y: x\sim y\} \vert$. A $d$-regular graph is a graph where every vertex has the same degree $d$. Let $A:=(a_{ij})$ be its adjacency matrix with $a_{xy}=1$ if $x \sim y$, otherwise $a_{xy}=0$. We denote the normalized adjacency matrix as
\[
P:=D^{-1}A,
\]
with $D$ the diagonal matrix of vertex degree: $D_{xx}=d_x$. We list the eigenvalues (counting multiplicity) of $P$ as
\[
1=\mu_1\ge\mu_2\ge\cdots\ge\mu_n\ge-1,
\]
here $n:=|V|$ is the number of vertices.

For a symmetric nonnegative matrix with spectral radius $\rho$
and eigenvalues ordered as in \eqref{eq:lambda}, write
\begin{equation}\label{eq:spectrum_gaps}
\gamma_j=\rho-\lambda_j,\qquad\delta_j=\rho+\lambda_{n-j+1},\qquad (1\le j \le n).
\end{equation}
Here $\delta_j$ is the $j$-th lower spectral gap and
$\gamma_j$ is the $j$-th upper spectral gap.
For a normalized adjacency matrix, these become
$1+\mu_{n-j+1}$ and $1-\mu_j$, respectively.

For a subset $B\subseteq V$,  its outer vertex boundary $\partial(B)$ is
\[
\partial_{out}(B):=
  \{x\in V\setminus B| \exists y\in B,\text{ such that }x\sim y\}.
\]

For a $d$-regular graph, the outer vertex isoperimetric constant $h_{out}(G)$ is defined as
\[
h_{out}(G):=\min_{\substack{\emptyset\neq B\subseteq V,\\|B|\le |V|/2}}\frac{|\partial_{out}(B)|}{|B|}.
\]

\begin{definition}
    Let $G=(V,E)$ be a finite simple graph with adjacency matrix $A$.
    We call $G$ \emph{even-walk-regular} if, for every even integer $\ell\ge 0$,
    $(A^\ell)_{xx}$ is independent of $x\in V$. If this holds for every
    integer $\ell\ge 0$, we call $G$ \emph{walk-regular}.
\end{definition}
Proposition~\ref{cg-prop:fiveEWR} in Appendix~\ref{app:five-classes} shows that the five classes of graphs which were investigated by Saha \cite{Saha} constitute a subclass of even-walk-regular graphs.

\subsection{Matrix}

We recall Perron--Frobenius theorem used below~\cite[Section~2.2]{BH12}.

A matrix $M=(M_{xy})\in\R^{n\times n}$ is \emph{nonnegative}
if $M_{xy}\ge0$ for all $x,y$.
Its directed support graph has vertex set $\{1,\ldots,n\}$
and an edge from $x$ to $y$ whenever $M_{xy}>0$.
The matrix is \emph{irreducible} if this directed graph is
strongly connected, that is, every vertex can be reached from
every other vertex by a directed path.

Let $M\in\mathbb{R}^{n\times n}$ and $\norm{M}$ denote the operator norm, that is
\[
\norm{M}=\sup_{x\in\mathbb{R}^n\setminus\{0\}}\frac{\norm{Mx}}{\norm{x}}.
\]

\begin{theorem}[Perron--Frobenius]\label{thm:perron-frobenius}
Let $M\in\R^{n\times n}$ be nonnegative, and let
\[
 \rho(M)=\max\{|\lambda|:\lambda\in\spec(M)\}
\]
be its spectral radius.
\begin{enumerate}
\item $\rho(M)$ is an eigenvalue of $M$, with a nonzero
eigenvector $h$ whose coordinates are all nonnegative.
\item If $M$ is irreducible and $M\ne0$, then $\rho(M)>0$
is algebraically simple and has an eigenvector $h$ whose
coordinates are all strictly positive.
\item If $Mh=rh$ for a strictly positive vector $h$, then
$r=\rho(M)$.
\end{enumerate}
If $M$ is symmetric, its largest eigenvalue satisfies
$\lambda_1(M)=\rho(M)=\norm M$.
\end{theorem}

The third assertion follows by putting $D_h=\diag(h_x)$.
The matrix $D_h^{-1}MD_h$ is nonnegative and has constant row
sum $r$. Thus $r$ is an eigenvalue and its spectral radius
is at most $r$, proving $\rho(M)=r$.

The Courant--Fischer theorem characterizes eigenvalues
through the Rayleigh quotient~\cite[Theorem~2.4.1]{BH12}.

\begin{theorem}[Courant--Fischer]\label{thm:courant-fischer}
Let $M$ be a self-adjoint operator on a $d$-dimensional real
inner product space $\mathcal H$, with eigenvalues
$\lambda_1(M)\ge\cdots\ge\lambda_d(M)$.
For every $1\le k\le d$,
\begin{align}\label{eq:courant-fischer}
 \lambda_k(M)
 &=\max_{\substack{L\subseteq\mathcal H\\\dim L=k}}
   \min_{\substack{F\in L\\\norm F=1}}\ip{F}{MF}\\
 &=\min_{\substack{L\subseteq\mathcal H\\\dim L=d-k+1}}
   \max_{\substack{F\in L\\\norm F=1}}\ip{F}{MF}.\notag
\end{align}
\end{theorem}

Ky Fan's  principle~\cite{Fan} characterizes the sum of the \(m\)
largest eigenvalues of a real symmetric matrix as a maximum over
rank-\(m\) orthogonal projections. Overton and Womersley give an
equivalent formulation in terms of positive semidefinite contractions
of trace \(m\) \cite[Eq.~(4)]{OvertonWomersley1992}.

\begin{theorem}[Ky Fan's principle]\label{thm:ky-fan}
Let \(N\in\mathbb{R}^{n\times n}\) be symmetric, with eigenvalues
\(\lambda_1\ge\cdots\ge\lambda_n\). For every integer \(1\le m\le n\)
and every symmetric matrix \(K\) satisfying \(0\preceq K\preceq I\)
and \(\operatorname{tr}(K)=m\), we have
\[
    \operatorname{tr}(NK)
    \le \sum_{j=1}^{m}\lambda_j.
\]
The upper bound is attained by projecting onto the span of eigenvectors
corresponding to the \(m\) largest eigenvalues.
\end{theorem}

\subsection{Topological theorem} We first recall a topological theorem due to Dold. Let $\Gamma$ be a group with identity $e$. A (left) action of $\Gamma$ on
a space $X$ assigns to each $g\in\Gamma$ a homeomorphism $x\mapsto gx$ of $X$
such that $ex=x$ and $(gh)x=g(hx)$ for all $g,h\in\Gamma$ and $x\in X$.
The action is \emph{free} if $gx\ne x$ for every $g\ne e$ and every $x\in X$.
If $\Gamma$ also acts on a space $Y$, a map $f\colon X\to Y$ is
\emph{$\Gamma$-equivariant} if $f(gx)=g f(x)$ for all $g\in\Gamma$ and $x\in X$.

\begin{theorem}[\cite{Dold}]\label{ed-thm:cyclic-borsuk-ulam}
Let $\Gamma$ be a nontrivial finite group acting freely on spaces $X$ and $Y$,
and let $n\ge 0$ be an integer. If $X$ is $n$-connected and $Y$ is a
paracompact Hausdorff space of covering dimension at most $n$, then there
is no continuous $\Gamma$-equivariant map $X\to Y$.
\end{theorem}


\section{Comparisons of lower and upper spectral gaps}\label{ed-sec:top}
In this section, we first consider the eigenvalue comparison for general symmetric matrices with each entry nonnegative, then apply it to the case of graphs.
Throughout this section, $N$ is a symmetric nonnegative matrix with
eigenvalues and operator norm as in \eqref{eq:lambda}, and its spectrum gaps $\{\gamma_j\}$ and $\{\delta_j\}$ as (\ref{eq:spectrum_gaps}).

\begin{theorem}[Lower and upper spectral gap comparisons]\label{ed-thm:main-matrix}
Let \(N\in\mathbb{R}^{n\times n}\) be a symmetric matrix with
nonnegative entries, with eigenvalues \(\lambda_1,\ldots,\lambda_n\)
and spectral radius \(\rho\) as in \eqref{eq:lambda}.
For every \(1\le m\le n\), we have
\begin{equation}\label{ed-eq:gap-scaled-matrix}
    \delta_m \ge \frac14\,\gamma_m
\end{equation}
and
\begin{equation}\label{eq:prefix-gaps-matrix}
    \sum_{j=1}^m \delta_j
    \ge \frac12\sum_{j=1}^m \gamma_j.
\end{equation}
\end{theorem}

We prove these comparisons of lower and upper spectral gaps using
a matrix inequality.

\subsection{Matrix lemma}\label{cg-sec:master}

Let $G$ be a complex Hermitian positive semidefinite matrix. We write
$G^*$ for its conjugate transpose and $\overline{G}$ for its complex
conjugate. Put
$r_x=\sqrt{G_{xx}}$ and $S=\{x:r_x>0\}$, and define two matrices
$K_0,K_2$ by
\[
(K_0)_{xy}=\frac{|G_{xy}|^2}{r_xr_y},\qquad
 (K_2)_{xy}=\frac{G_{xy}^{\,2}}{r_xr_y}
 \quad (x,y\in S),
\]
and set their other entries to zero.

For real $G$, we have $K_0=K_2$. The construction of $K_0$ in this case
appears, with trace-one scaling, in
\cite[Lemma~IV.4]{BH}.

\begin{lemma}\label{cg-lem:master}
Let $G=G^*\succeq0$ be a complex Hermitian positive semidefinite matrix.
Then the matrices $K_0,K_2$ defined above satisfy
\begin{equation}\label{cg-eq:structural}
 K_0,K_2\succeq0,\qquad
 \Tr K_0=\Tr K_2=\Tr G=\norm r^2.
\end{equation}
For every real symmetric nonnegative matrix $N$ of the same size, it holds that
\begin{equation}\label{cg-eq:master}
 2|\Tr(NG)|^2\le
 (r^{\mathsf T}Nr)\bigl[\Tr(NK_0)+\Tr(NK_2)\bigr].
\end{equation}
Moreover, if $G\preceq I$, then $K_0,K_2\preceq I$.
\end{lemma}

\begin{proof}
By positivity of $G$, we have $|G_{xy}|^2\le G_{xx}G_{yy}$, so every row and column outside $S$ vanishes in $K_0$ and $K_2$. Choose a factorization
$G=\sum_\ell a_\ell a_\ell^*$ with $a_\ell\in\C^n$, and define diagonal matrices $D_\ell$ by
\[
 (D_\ell)_{xx}=\frac{a_\ell(x)}{r_x}\quad(x\in S),\qquad
 (D_\ell)_{xx}=0\quad(x\notin S).
\]
If $I_S$ is the coordinate projection onto $S$, entrywise calculation gives
\[
\sum_\ell D_\ell D_\ell^*=I_S,\qquad
 K_0=\sum_\ell D_\ell\overline G D_\ell^*,\qquad
 K_2=\sum_\ell D_\ell G D_\ell^*.
\]
Since both $G$ and $\overline G$ are positive semidefinite, the above identities prove positivity of $K_0$ and $K_2$. Also,
$K_j\preceq I_S\preceq I$ is obtained by $G\preceq I$. Their diagonal entries equal $G_{xx}$, which proves the trace identity in \eqref{cg-eq:structural}.

With every real symmetric nonnegative matrix $N$ of the same size, put $\alpha=r^{\mathsf T}Nr$ and $T=\Tr(NG)$. For $x,y\in S$,
define
\[
 w_{xy}=N_{xy}r_xr_y\ge0,\qquad
 c_{xy}=\frac{\operatorname{Re}G_{xy}}{r_xr_y}\in\R.
\]
The trace of a product of Hermitian matrices is real. Since $N$ is real
symmetric, conjugate terms cancel in the sums, and
\[
 \alpha=\sum_{x,y\in S}w_{xy},\qquad
 T=\sum_{x,y\in S}w_{xy}c_{xy}.
\]
The scalar identity $|z|^2+\operatorname{Re}(z^2)=2(\operatorname{Re}z)^2$
gives
\[
\frac{\Tr(NK_0)+\Tr(NK_2)}{2}
 =\sum_{x,y\in S}w_{xy}c_{xy}^2.
\]
By Cauchy--Schwarz with the nonnegative weights $w_{xy}$,
\begin{align*}
 T^2
 &=\left(\sum_{x,y\in S}w_{xy}c_{xy}\right)^2\le
   \left(\sum_{x,y\in S}w_{xy}\right)
   \left(\sum_{x,y\in S}w_{xy}c_{xy}^2\right)\\
  &=\frac{\alpha}{2}\bigl[\Tr(NK_0)+\Tr(NK_2)\bigr],
\end{align*}
which leads to \eqref{cg-eq:master}. When $\alpha>0$, equality holds precisely
when $c_{xy}$ is constant on the pairs with $x,y\in S$ and
$N_{xy}>0$. When $\alpha=0$, all weights vanish, so
$T=\Tr(NK_0)=\Tr(NK_2)=0$ and equality is automatic.
\end{proof}

\subsection{Comparison of individual spectral gaps}
On $\C^n$, we use $\ip zw=\sum_{i=1}^n\overline{z_i}w_i$ as the inner product of two vectors $z=(z_1,z_2,\cdots,z_n),w=(w_1,w_2,\cdots,w_n)\in \mathbb{C}^n$. And we define a function $\phi\colon\mathbb{C} \to \mathbb{C}$ as
\[
\phi(w)=\begin{cases}\frac{w^2}{|w|},&w\ne0,\\0,&w=0,\end{cases}
\]
and a map $\Phi\colon\mathbb{C}^n \to \mathbb{C}^n$ by
\[
(\Phi(z))_i=\phi(z_i), \qquad \text{for} \quad i=1,2,\cdots,n\quad \text{and} \quad z=(z_1,z_2,\cdots,z_n) \in \mathbb{C}^n.
\]
The map $\Phi$ is continuous, including at zero, and it satisfies
\[
\norm{\Phi(z)}=\norm z,\qquad
 \Phi(e^{it}z)=e^{2it}\Phi(z)\quad(t\in\R).
\]
Moreover, we have the following observation regarding the map $\Phi$.
\begin{corollary}\label{ed-cor:rayleigh}
Let $N$ be a real symmetric nonnegative matrix with
$\rho>0$. Every unit vector $z\in\C^n$ satisfies
\[
\ip{\Phi(z)}{N\Phi(z)}\ge
 \frac{2\ip z{Nz}^{\,2}}{\rho}-\rho.
\]
\end{corollary}
\begin{proof}
Apply Lemma~\ref{cg-lem:master} with $G=zz^*$. This gives
\[
r=|z|,\qquad K_0=rr^{\mathsf T},\qquad
 K_2=\Phi(z)\Phi(z)^*.
\]
Writing $T=\ip z{Nz}$, $B=\ip{\Phi(z)}{N\Phi(z)}$, and
$\alpha=r^{\mathsf T}Nr$, we have $0\le\alpha\le\rho$ because $r$
is a nonnegative unit vector. Equation~\eqref{cg-eq:master} gives
\[
2T^2\le\alpha(\alpha+B).
\]
If $\alpha=0$, then $T=B=0$ and the claim is immediate. Otherwise we have
\[
 B\ge\frac{2T^2}{\alpha}-\alpha
 \ge\frac{2T^2}{\rho}-\rho.
\]
This proves the corollary.
\end{proof}
We now apply the sphere case of Dold's theorem (Theorem~\ref{ed-thm:cyclic-borsuk-ulam}).

\begin{lemma}\label{ed-lem:orthogonal}
Let $H,W\subseteq\C^n$ be complex linear subspaces with
\[
 \dim_{\C}H=m,\qquad \dim_{\C}W=r<m.
\]
There exists $z\in H$ such that $\norm{z}=1$ and $\Phi(z)\perp W$.
\end{lemma}

\begin{proof}
If $r=0$, any unit vector in $H$ suffices. Now assume $r\ge1$, and denote by
$S(H)$ and $S(W)$ the unit spheres in $H$ and $W$, respectively. Let $\Pi_W$ be the complex orthogonal projection onto $W$. Consider the continuous map
\[
 T\colon S(H)\to W,\qquad T(z)=\Pi_W\Phi(z).
\]
Let $\zeta=e^{2\pi i/3}$. Since $\Pi_W$ is complex linear, we have
\begin{equation}\label{ed-eq:equivariance}
 T(\zeta z)=\zeta^2 T(z).
\end{equation}

Suppose, for a contradiction, that $T$ has no zero, we can define a continuous map $\widehat T$ by normalizing the vector $T(z)$ as follows
\[
 \widehat T\colon S(H)\to S(W),
 \qquad \widehat T(z)=\frac{T(z)}{\norm{T(z)}}.
\]
Let $\Gamma=\langle\zeta\rangle$ act on $S(H)$ by $\zeta\cdot z=\zeta z$
and on $S(W)$ by $\zeta\cdot w=\zeta^2w$. Each nonidentity element acts
on each sphere by a scalar different from $1$, so both actions are free.
Equation~\eqref{ed-eq:equivariance} makes $\widehat T$ $\Gamma$-equivariant.
However, the dimensions of the unit spheres are
\[
 S(H)\cong S^{2m-1},\qquad S(W)\cong S^{2r-1},
 \qquad 2m-1>2r-1,
\]
contradicting the sphere case of Theorem~\ref{ed-thm:cyclic-borsuk-ulam}. Hence $T$ has a zero.
\end{proof}

We first compare individual lower and upper spectral gaps in Theorem~\ref{ed-thm:main}.
\begin{proof}[Proof of (\ref{ed-eq:gap-scaled-matrix}) in Theorem~\ref{ed-thm:main-matrix}]
We prove the stronger estimate
\begin{equation}\label{ed-eq:root-scaled}
 \lambda_{n-m+1}\ge-\sqrt{\frac{\rho(\rho+\lambda_m)}{2}}.
\end{equation}
If $\rho=0$, then $N=0$ and all conclusions are immediate. Assume
$\rho>0$ and fix $m$. If $\lambda_{n-m+1}\ge0$, then
\eqref{ed-eq:root-scaled} is immediate because $\rho+\lambda_m\ge0$.
Otherwise write $\lambda_{n-m+1}=-a<0$.

Choose a real orthonormal eigenbasis $v_1,\ldots,v_n$ for $N$ and regard
it as a complex orthonormal basis of $\C^n$. Set
\[
 H=\Span_{\C}\{v_{n-m+1},\ldots,v_n\},
 \qquad
 W=\Span_{\C}\{v_1,\ldots,v_{m-1}\}.
\]
Their complex dimensions are $m$ and $m-1$, respectively. Every unit
vector $z\in H$ satisfies
\begin{equation}\label{ed-eq:smallest-rayleigh}
 \ip{z}{Nz}\le-a.
\end{equation}
By Lemma~\ref{ed-lem:orthogonal}, choose such a $z$ with
\[
 f:=\Phi(z)\perp W.
\]
Then $\norm{f}=\norm{z}=1$. Corollary~\ref{ed-cor:rayleigh} and the inequality
\eqref{ed-eq:smallest-rayleigh} imply
\[
 \ip{f}{Nf}\ge\frac{2\ip{z}{Nz}^{\,2}}{\rho}-\rho
 \ge\frac{2a^2}{\rho}-\rho.
\]
On the other hand, expansion in the eigenbasis yields
\[
 \ip{f}{Nf}
 =\sum_{j=m}^n\lambda_j|\ip{v_j}{f}|^2
 \le\lambda_m\sum_{j=m}^n|\ip{v_j}{f}|^2
 =\lambda_m.
\]
Thus $2a^2\le\rho(\rho+\lambda_m)$, proving
\eqref{ed-eq:root-scaled}.
Since $0\le\gamma_m\le2\rho$, the root bound gives
\[
 \delta_m\ge\rho-\sqrt{\rho\left(\rho-\frac{\gamma_m}{2}\right)}
 =\frac{\gamma_m}{2\left(1+\sqrt{1-\frac{\gamma_m}{2\rho}}\right)}
 \ge\frac{\gamma_m}{4},
\]
which is \eqref{ed-eq:gap-scaled}.
\end{proof}

\subsection{Sums over consecutive eigenvalues}\label{sw-sec:windows}
\begin{proof}[Proof of (\ref{eq:prefix-gaps-matrix}) in Theorem~\ref{ed-thm:main-matrix}]
If $\rho=0$, the claim is immediate. Assume $\rho>0$.
Let $Q$ be a real rank-$m$ orthogonal projection. Applying
Lemma~\ref{cg-lem:master} with $G=Q$, we obtain
\[
0\preceq K_i\preceq I \qquad  \text{and }\qquad  \Tr K_i=m \qquad\text{ for }i=0,2.
\]
Ky Fan's principle gives
\[
\Tr(NK_i)\le\sum_{j=1}^m\lambda_j.
\]
In an orthonormal eigenbasis of
$N$, the diagonal entries of $K_i$ lie in $[0,1]$ and sum to $m$.
Since $Q$ is real, $K_0=K_2$.
Writing $r_x=\sqrt{Q_{xx}}$ and
$S=\{x:r_x>0\}$, we have
\[
 \Tr(N(K_0+K_2))
 =2\sum_{x,y\in S}N_{xy}\frac{Q_{xy}^{\,2}}{r_xr_y}\ge0.
\]
Since $r^{\mathsf T}Nr\le m\rho$, it follows that
\begin{equation}\label{eq:kyfan}
 |\Tr(NQ)|^2
 \le\tfrac12(r^{\mathsf T}Nr)\Tr(N(K_0+K_2))
 \le m\rho\sum_{j=1}^m\lambda_j.
\end{equation}
Taking $Q$ to project onto the bottom $m$ eigenvectors gives
\begin{equation}\label{eq:prefix-eigenvalues}
 \left(\sum_{j=n-m+1}^n\lambda_j\right)^2
 \le m\rho\sum_{j=1}^m\lambda_j.
\end{equation}
Rewriting \eqref{eq:prefix-eigenvalues} in terms of the gaps gives
\[
 \left(\sum_{j=1}^m\delta_j-m\rho\right)^2
 \le m\rho\left(m\rho-\sum_{j=1}^m\gamma_j\right).
\]
The lower square-root bound yields
\[
 \sum_{j=1}^m\delta_j
 \ge m\rho\left(1-\sqrt{1-
        \frac{\sum_{j=1}^m\gamma_j}{m\rho}}\right)
 \ge\frac12\sum_{j=1}^m\gamma_j,
\]
where $0\le\sum_{j=1}^m\gamma_j\le m\rho$ follows from
\eqref{eq:kyfan}, and $1-\sqrt{1-t}\ge t/2$ for $0\le t\le1$.
This proves \eqref{eq:prefix-gaps}.
\end{proof}

\begin{corollary}[Sums over consecutive eigenvalues]\label{cor:windows}
Let $k\ge0$ and $q\ge1$ be integers with $k+q\le n$. Then
\begin{equation}\label{sw-eq:combinedwindow}
 \sum_{j=k+1}^{k+q}\delta_j
 \ge\frac{q}{2(q+\min\{k,q\})}
      \sum_{j=k+1}^{k+q}\gamma_j.
\end{equation}
\end{corollary}
\begin{proof}
Put $m=k+q$. Since the $\delta_j$ are nondecreasing,
\[
 \sum_{j=k+1}^{m}\delta_j
 \ge\frac{q}{m}\sum_{j=1}^m\delta_j
 \ge\frac{q}{2m}\sum_{j=1}^m\gamma_j
 \ge\frac{q}{2m}\sum_{j=k+1}^m\gamma_j,
\]
where the middle inequality is \eqref{eq:prefix-gaps-matrix}.
Also \eqref{ed-eq:gap-scaled-matrix} gives
$\sum_{j=k+1}^m\delta_j\ge
\frac14\sum_{j=k+1}^m\gamma_j$.
The larger of the coefficients $q/(2m)$ and $1/4$ is
$q/(2(q+\min\{k,q\}))$.
\end{proof}

\begin{corollary}\label{cor:cross-rank}
For $1\le k\le\ell\le n$,
\begin{equation}\label{eq:cross-rank}
 \delta_\ell\ge
 \max\left\{\frac14,\frac{\ell-k+1}{2\ell}\right\}\gamma_k.
\end{equation}
\end{corollary}
\begin{proof}
Monotonicity and \eqref{eq:prefix-gaps-matrix} give
\[
 \ell\delta_\ell\ge\sum_{j=1}^\ell\delta_j
 \ge\frac12\sum_{j=1}^\ell\gamma_j
 \ge\frac{\ell-k+1}{2}\gamma_k.
\]
Also \eqref{ed-eq:gap-scaled-matrix} and $\gamma_\ell\ge\gamma_k$ give
$\delta_\ell\ge\gamma_k/4$.
\end{proof}
\subsection{Counting eigenvalues at the lower and upper ends}\label{sec:threshold-counting}
Following \cite{BH}, for a real symmetric matrix $N$ with eigenvalues
$\lambda_1\ge\cdots\ge\lambda_n$ and $a\in\R$, write
\[
 \operatorname{rank}_{\le a}(N):=\#\{i:\lambda_i\le a\},
 \qquad
 \operatorname{rank}_{\ge a}(N):=\#\{i:\lambda_i\ge a\},
\]
where eigenvalues are counted with multiplicity.
The partial-sum estimate \eqref{eq:prefix-eigenvalues}
gives the following strengthening of \cite[Theorem~4.1]{BH}.

\begin{corollary}[Large lower rank implies large upper rank]
\label{cor:threshold-counting}
Let $N\in\R^{n\times n}$ be a symmetric matrix with nonnegative
entries and $\norm N\le1$. If
\[
 \operatorname{rank}_{\le-\tau}(N)\ge t
\]
for some $\tau\ge0$ and some positive integer $t$, then, for every
$\sigma\in(0,1)$,
\begin{equation}\label{eq:threshold-counting}
 \operatorname{rank}_{\ge\frac{\tau^2-\sigma}{1-\sigma}}(N)
 \ge\sigma t.
\end{equation}
\end{corollary}

\begin{proof}
If $N=0$, the hypothesis forces $\tau=0$ and the conclusion is
immediate. Assume $N\ne0$ and put $\rho=\norm N$.
Then $0<\rho\le1$, $0\le\tau\le1$, and each of the  $t$ lowest
eigenvalues is at most $-\tau$. Thus \eqref{eq:prefix-eigenvalues}
with $m=t$ gives
\[
 \sum_{j=1}^t\lambda_j
 \ge\frac{1}{t\rho}
      \left(\sum_{j=n-t+1}^n\lambda_j\right)^2
 \ge\frac{t\tau^2}{\rho}
 \ge t\tau^2.
\]
If $\tau=1$, this forces $\lambda_1=\cdots=\lambda_t=1$,
and the conclusion follows.
For $0\le\tau<1$, set
\[
 b:=\frac{\tau^2-\sigma}{1-\sigma}<1,
 \qquad s:=\min\{t,\operatorname{rank}_{\ge b}(N)\}.
\]
Among the largest $t$ eigenvalues, exactly $s$ are at least $b$.
Each eigenvalue is at most $1$, so
\[
 t\tau^2\le\sum_{j=1}^t\lambda_j
 \le s+(t-s)b.
\]
Rearranging and using the definition of $b$ gives
\[
 \operatorname{rank}_{\ge b}(N)
 \ge s\ge\frac{\tau^2-b}{1-b}\,t=\sigma t,
\]
as required.
\end{proof}

Buhai, Hua, Steurer and V{\'a}ri-Kakas~\cite[Theorem~4.1]{BH}
obtain the lower bound $\sigma^2t$ under the same hypotheses. Thus \eqref{eq:threshold-counting}
improves the coefficient from $\sigma^2$ to $\sigma$.

\section{Two-sided spectral gap comparisons and applications}
\label{rf-sec:angular}
\label{sec:applications}

\subsection{Walk-regular graphs}
We compare eigenvalues at opposite ends of the spectrum using a
vector-valued eigenfunction whose length agrees with a positive
Perron eigenvector.

For an integer $r\ge1$, put $q=2^r$ and set
\[
 \begin{gathered}
  X=\begin{pmatrix}0&1\\1&0\end{pmatrix},\qquad
  Z=\begin{pmatrix}1&0\\0&-1\end{pmatrix},\\[0.4ex]
  C_i=Z^{\otimes(i-1)}\otimes X\otimes I_2^{\otimes(r-i)}
  \quad(1\le i\le r).
 \end{gathered}
\]
These real symmetric matrices satisfy
\begin{equation}\label{eq:matrix-relations}
 C_iC_j+C_jC_i=2\delta_{ij}I_q.
\end{equation}
Indeed, $C_i^2=I_q$, and for $i<j$ anticommutation follows from
$XZ=-ZX$ at tensor position $i$.

\begin{theorem}\label{thm:perron-angular}
Let $N\ne0$ be an irreducible symmetric nonnegative matrix indexed by $V$, with eigenvalues $\{\lambda_i\}$ and the spectrum radius $\rho$ as in \eqref{eq:lambda}, and let
$h\in(0,\infty)^n$ satisfy $Nh=\rho h$. Suppose $0<a\le1$ and a map $u\colon V\to\R^r$ satisfies
\begin{equation}\label{eq:perron-map}
 Nu=-a\rho u,\qquad |u(x)|=h_x\quad(x\in V),
\end{equation}
where $N$ acts on each coordinate of $u$. Then, for $1\le j\le n$,
\begin{equation}\label{eq:perron-angle}
 \left|\arccos\frac{\lambda_j}{\rho}
       +\arccos\frac{\lambda_{n-j+1}}{\rho}-\pi\right|
 \le\arccos a.
\end{equation}
One also has
\begin{equation}\label{eq:perron-roots}
 |\sqrt{\gamma_j}-\sqrt{\delta_j}|
 \le\sqrt{\rho(1-a)},\qquad
 |\delta_j-\gamma_j|\le\rho\sqrt{2(1-a)}.
\end{equation}
\end{theorem}

\begin{proof}
Put $D_h=\diag(h_x)$ and
\[
 P=\rho^{-1}D_h^{-1}ND_h,\qquad
 \pi_x=\frac{h_x^2}{\sum_yh_y^2},\qquad
 v(x)=\frac{u(x)}{h_x}.
\]
Then $P$ is stochastic and reversible with respect to $\pi$.
Its eigenvalues are $\mu_j=\lambda_j/\rho$, and
$|v(x)|=1$, $Pv=-av$.

Let $\mathcal H=L^2(V,\pi;\R^q)$, let $\mathcal P=P\otimes I_q$,
and define
\[
 U_x=\sum_{i=1}^rv_i(x)C_i,\qquad
 (\mathcal UF)(x)=U_xF(x).
\]
The relations in \eqref{eq:matrix-relations} and $Pv=-av$ give
\begin{equation}\label{rf-eq:Uproperties}
 U_x^{\mathsf T}=U_x,\qquad U_x^2=I_q,\qquad
 \sum_yP_{xy}U_y=-aU_x.
\end{equation}
Thus $\mathcal U$ is an orthogonal involution on $\mathcal H$.

For a unit vector $F\in\mathcal H$, write
\[
 t=\ip{F}{\mathcal PF},\qquad
 s=\ip{\mathcal UF}{\mathcal P\mathcal UF}.
\]
Put $\omega_{xy}=\pi_xP_{xy}$. Both marginals of $\omega$ equal $\pi$,
so the following vectors are unit vectors in
$L^2(V\times V,\omega;\R^q)$:
\[
 A(x,y)=-F(x),\qquad
 B(x,y)=U_yU_xF(x),\qquad
 C(x,y)=F(y).
\]
Using \eqref{rf-eq:Uproperties}, their inner products are
\[
 \ip{A}{B}_{\omega}=a,\qquad
 \ip{A}{C}_{\omega}=-t,\qquad
 \ip{B}{C}_{\omega}=s.
\]
The reverse triangle inequality for angles between unit vectors gives
\begin{equation}\label{rf-eq:pointwise-angle}
 \left|\arccos s-\arccos(-t)\right|\le\arccos a.
\end{equation}

For a self-adjoint contraction $M$,
Theorem~\ref{thm:courant-fischer} and the monotonicity of
$\arccos$ give
\[
 \arccos\lambda_k(M)
 =\min_{\substack{L\subseteq\mathcal H\\\dim L=k}}
   \max_{\substack{F\in L\\\norm{F}=1}}
   \arccos\ip{F}{MF}.
\]
Apply this formula to $\mathcal U\mathcal P\mathcal U$ and
$-\mathcal P$ in \eqref{rf-eq:pointwise-angle}. Their eigenvalues are,
respectively, $\mu_1,\ldots,\mu_n$ and $-\mu_n,\ldots,-\mu_1$,
each repeated $q$ times. Taking $k=jq$ yields
\[
 \left|\arccos\mu_j-\arccos(-\mu_{n-j+1})\right|\le\arccos a,
\]
which proves \eqref{eq:perron-angle}.

Finally, put
$\alpha=\arccos\mu_j$, $\beta=\arccos(-\mu_{n-j+1})$, and
$\delta=\arccos a$. Then $|\alpha-\beta|\le\delta$.
For $b,c\in[0,\pi/2]$, one has
$|\sin b-\sin c|\le\sin|b-c|$. Hence,
\begin{align*}
 |\sqrt{\gamma_j}-\sqrt{\delta_j}|
 &=\sqrt{2\rho}\left|\sin\frac{\alpha}{2}
                         -\sin\frac{\beta}{2}\right|\\
 &\le\sqrt{2\rho}\sin\frac{\delta}{2}
   =\sqrt{\rho(1-a)},\\
 |\delta_j-\gamma_j|
 &=\rho|\cos\alpha-\cos\beta|\\
 &\le2\rho\sin\frac{\delta}{2}
   =\rho\sqrt{2(1-a)}.
\end{align*}
This proves \eqref{eq:perron-roots}.
\end{proof}

\begin{remark}
    Spherical eigenmaps have been used to derive universal eigenvalue
inequalities for Schr\"odinger operators on Riemannian manifolds;
see El Soufi, Harrell and Ilias~\cite[Section~4]{MR2471921}.
Here we use unit-length eigenmaps to compare eigenvalues at opposite
ends of the spectrum.
\end{remark}

Let $G$ be a finite simple $d$-regular graph on $n\ge2$ vertices,
with $d\ge1$, and let $A$ denote its adjacency matrix.
Write $1=\mu_1\ge\cdots\ge\mu_n\ge-1$ for the eigenvalues of
its normalized adjacency matrix $P=A/d$.
The $j$-th upper and lower spectral gaps of $P$ are, respectively,
\[
 \gamma_j=1-\mu_j,\qquad
 \delta_j=1+\mu_{n-j+1}\qquad(1\le j\le n).
\]

We recall the following classical property of walk-regular
graphs and include a proof for completeness; see
\cite{GM,CamaraVanDamKoolenPark2013}.

\begin{proposition}\label{rf-prop:walk-projections}
If $G$ is walk-regular, every spectral projection of $P$ has constant
diagonal.
\end{proposition}
\begin{proof}
If $\Theta$ is the set of distinct eigenvalues of $P$, spectral
interpolation gives
\[
 E_\theta=\prod_{\substack{\beta\in\Theta\\\beta\ne\theta}}
                  \frac{P-\beta I}{\theta-\beta}.
\]
Walk-regularity makes every polynomial in $P$ have constant
diagonal. Since $\Tr E_\theta=\operatorname{rank}E_\theta$,
the common diagonal entry is $\operatorname{rank}(E_\theta)/n$.
\end{proof}

For any spectral projection $E_\theta$ of a walk-regular graph,
let $r=\operatorname{rank}(E_\theta)$ and choose an orthonormal basis
$f_1,\ldots,f_r$ of the $\theta$-eigenspace of $P$. Then
\begin{equation}\label{rf-eq:mapfromprojection}
 v(x)=\sqrt{\frac nr}\,(f_1(x),\ldots,f_r(x))
\end{equation}
has unit length and satisfies $Pv=\theta v$.

\begin{proof}[Proof of Theorem~\ref{rf-thm:rootgaps}]
A nonzero simple graph has trace-zero adjacency and hence
$\mu_n<0$. For a walk-regular graph, Proposition~\ref{rf-prop:walk-projections}
and \eqref{rf-eq:mapfromprojection} provide a unit eigenmap for
$\mu_n$. Theorem~\ref{thm:perron-angular}
with $N=P$, $h=\mathbf1$, and $a=-\mu_n$ proves the theorem.
\end{proof}

\begin{remark}
    There are many natural examples of the hypothesis of Theorem~\ref{thm:perron-angular}. Every
vertex-transitive graph is walk-regular, since automorphisms
conjugate its adjacency matrix and its powers. This includes
Cayley graphs, complete graphs, and cycles.
In a connected bipartite regular graph, the signs of the two parts
give a scalar eigenmap with $a=1$. More generally, every irreducible symmetric
nonnegative matrix whose support is bipartite has the Perron
eigenmap $u(x)=\sigma_xh_x$, again with $a=1$.
\end{remark}

\subsection{Even-walk-regular graphs}

\begin{proposition}\label{rf-prop:EWR-projections}
Let $G$ be a finite connected even-walk-regular graph. Every spectral projection
of $P^2$ has constant diagonal. If $-\mu_n\notin\spec(P)$, the
smallest-eigenvalue projection $E_{\mu_n}$ has constant diagonal and all
conclusions of Theorem~\ref{rf-thm:rootgaps} hold.
\end{proposition}
\begin{proof}
Each spectral projection of $P^2$ is a polynomial in $P^2$, by
interpolation on its distinct eigenvalues. Even-walk-regularity makes the diagonal
of every such polynomial constant. Let $Q_{\mu_n^2}$ denote the
spectral projection of $P^2$ onto its $\mu_n^2$-eigenspace. Then
\[
 Q_{\mu_n^2}=E_{\mu_n}+E_{-\mu_n},
\]
where an absent eigenvalue contributes zero.
If $-\mu_n$ is absent, then $Q_{\mu_n^2}=E_{\mu_n}$, so
$E_{\mu_n}$ has constant diagonal.
Formula \eqref{rf-eq:mapfromprojection} supplies the eigenmap.
\end{proof}

The five classes of graphs considered by Saha~\cite{Saha} are
even-walk-regular, as shown in Proposition~\ref{cg-prop:fiveEWR}
in Appendix~\ref{app:five-classes}.

The balanced equitable partition associated with a nontrivial simple
eigenvalue of a connected walk-regular graph is classical; see Godsil
and McKay~\cite[Theorems~4.1 and~4.6]{GM} and
van Dam~\cite[Theorem~3.3]{MR1344560}.
The following lemma requires only the constant-diagonal property of
the rank-one spectral projection and additionally establishes the
comparison \eqref{rf-eq:linear} for the entire spectrum.

\begin{lemma}\label{rf-lem:linear}
Suppose that $G$ is a connected nonbipartite regular simple graph, the smallest eigenvalue $\mu_n$ is simple, and its projection has constant diagonal. Then there is a balanced equitable bipartition $V=V_+\sqcup V_-$ and an integer $\ell$ with
\[
1\le\ell<\frac{d}{2},\qquad \delta_1=\frac{2\ell}{d},
\]
such that every vertex has exactly $\ell$ neighbors in its own cell. For every $1\le j\le n$,
\begin{equation}\label{rf-eq:linear}
 |\delta_j-\gamma_j|\le\delta_1.
\end{equation}
Equivalently,
\begin{equation}\label{rf-eq:linearbracket}
 \max\{\delta_1,\gamma_j-\delta_1\}\le\delta_j\le\gamma_j+\delta_1,
 \qquad \text{and }\qquad \delta_j\ge\frac{\gamma_j}{2}.
\end{equation}
Moreover, if $n\equiv2\pmod4$, then the integer $\ell$ is even. If $n$ is odd, the smallest eigenvalue cannot be simple with a constant-diagonal spectral projection.
\end{lemma}
\begin{proof}
Let $f$ be a real eigenvector for $\mu_n$, normalized by $\sum_{x\in V}f(x)^2=1$. Since its rank-one projection has constant diagonal, \[
f(x)^2=\frac{1}{n},\qquad \text{for every }x\in V.
\]
Thus $\sigma:=\sqrt n f$ satisfies $\sigma(x)\in\{+1,-1\}$. Since $f\perp\mathbf{1}$, the sign classes of $\sigma$ are balanced, given by $\sum_{x\in V}\sigma(x)=\sqrt{n}\sum_{x\in V}f(x)=0$. Let $\ell_x$ be the number of neighbors of $x$ with the same sign, the eigenfunction equation gives
\[
 d\mu_n\sigma(x)=(2\ell_x-d)\sigma(x).
\]
Hence we derive that
\begin{equation}\label{eq:l_epsilon_1}
\ell_x=\frac{d(1+\mu_n)}{2}=\frac{d\delta_1}{2}
\end{equation}
is independent of vertices, i.e. the same $\ell$ at every vertex $x$. Since the graph is nonbipartite, we have $\ell>0$, while $\mu_n<0$ yields $\ell<\frac{d}{2}$.

Put $S=\diag(\sigma)$ and define $P_{\mathrm{in}}$ entrywise by
\[
 (P_{\mathrm{in}})_{xy}=
 \begin{cases}
  P_{xy},&\sigma(x)=\sigma(y),\\
  0,&\sigma(x)\ne\sigma(y).
 \end{cases}
\]
Conjugation by $S$ preserves same-sign entries and negates
opposite-sign entries, so
\begin{equation}\label{rf-eq:switch}
 SPS=-P+2P_{\mathrm{in}}.
\end{equation}
Since $S\mathbf1=\sigma$, $P\sigma=\mu_n\sigma$, and
$P\mathbf1=\mathbf1$, applying this identity to $\mathbf1$ gives
\[
 2P_{\mathrm{in}}\mathbf1
 =(1+\mu_n)\mathbf1=\delta_1\mathbf1.
\]
Since the symmetric nonnegative matrix $2P_{\mathrm{in}}$ has the positive
eigenvector $\mathbf1$ with eigenvalue $\delta_1$,
Theorem~\ref{thm:perron-frobenius} gives $\norm{2P_{\mathrm{in}}}=\delta_1$.
As $SPS$ and $P$ have the same eigenvalues, Weyl's inequality applied
to \eqref{rf-eq:switch} yields
\[
 |\mu_j+\mu_{n-j+1}|\le\delta_1,
\]
which is \eqref{rf-eq:linear}. Since $\delta_j\ge\delta_1$,
the bound also gives \eqref{rf-eq:linearbracket}; its last inequality
follows from $\max\{\delta_1,\gamma_j-\delta_1\}\ge\gamma_j/2$.

Finally, if $n/2$ is odd, the induced graph $G[V_+]$ is
$\ell$-regular of odd order. The handshaking identity gives
\[
\ell \cdot \frac{n}{2}=2|E(G[V_+])|.
\]
Thus $\ell$ is even when $n/2$ is odd. In addition, if $n$ is odd, a balanced sign partition cannot exist.
\end{proof}

\begin{proof}[Proof of Theorem~\ref{thm:ewr-gap}]
Set $F(s)=1-\sqrt{1-s/2}$ for $0\le s\le2$.
We prove the stronger bound
\[
 \delta_1\ge\min\left\{\frac2d,F(\gamma_2)\right\}.
\]
If $\delta_1\ge F(\gamma_2)$, this stronger bound is immediate.
Suppose instead that
\begin{equation}\label{eq:strict_inequality}
 \delta_1<F(\gamma_2).
\end{equation}
If $\mu_n$ were multiple, \eqref{ed-eq:root-scaled} with $m=2$ and $\rho=1$
would give $\delta_1=\delta_2\ge F(\gamma_2)$, a
contradiction. Thus $\mu_n$ is simple. Since $F(s)\le s$ on
$[0,2]$, \eqref{eq:strict_inequality} implies $-\mu_n>\mu_2$.
Nonbipartiteness gives $-\mu_n<1=\mu_1$, so
$-\mu_n\notin\spec(P)$. By
Proposition~\ref{rf-prop:EWR-projections}, the rank-one projection
$E_{\mu_n}$ has constant diagonal.
Lemma~\ref{rf-lem:linear} then gives a balanced equitable
bipartition $V=V_+\sqcup V_-$ in which every vertex has exactly
$\ell$ neighbors in its own cell, where
$1\le\ell<d/2$ and $\delta_1=2\ell/d\ge2/d$.
It also gives $|\delta_j-\gamma_j|\le\delta_1$ for every
$j$. This proves the stronger bound.

The strict case \eqref{eq:strict_inequality} is impossible if
$n$ is odd or no balanced equitable bipartition exists. In either
situation, $\delta_1\ge F(\gamma_2)$. If
$n\equiv2\pmod4$, Lemma~\ref{rf-lem:linear} makes $\ell$ even
in the strict case, hence $\ell\ge2$ and $\delta_1\ge4/d$.
In the other case $\delta_1\ge F(\gamma_2)$. Thus, when
$n\equiv2\pmod4$,
\[
 \delta_1\ge\min\left\{\frac4d,F(\gamma_2)\right\}.
\]
Finally, $F(s)\ge s/4$ for $0\le s\le2$, giving the first
inequality in \eqref{eq:ewr-gap}. Also $d\ge2$ and $0<\gamma_2\le2$, so
$2/d\ge\gamma_2/(2d)$ and
$F(\gamma_2)\ge\gamma_2/4\ge\gamma_2/(2d)$.
This proves the second inequality in \eqref{eq:ewr-gap}.
\end{proof}

\section{Sharpness of the comparisons}\label{ed-sec:sharp}

We examine the constants in the comparisons of individual spectral
gaps and their partial sums.

The next proposition identifies all pairs of indices for which the factor
$1/4$ is sharp, even among connected vertex-transitive graphs.

\begin{proposition}[Sharpness across indices]\label{ed-prop:sharp-gap}
Fix integers $k\ge2$ and $\ell\ge1$. Among finite connected nonbipartite
vertex-transitive simple graphs $G$ with $n=|V(G)|\ge\max\{k,\ell\}$,
\begin{equation}\label{ed-eq:sharp-gap}
 \inf_G\frac{\delta_\ell}{\gamma_k}=\frac14
 \quad\text{if and only if}\quad k\le\ell\le2k-2.
\end{equation}
If $\ell<k$, then $\inf_G \frac{\delta_\ell}{\gamma_k}=0$.
\end{proposition}

\begin{proof}\leavevmode
\begin{enumerate}
\renewcommand{\labelenumi}{\textup{(\roman{enumi})}}
\setlength{\itemsep}{0.6\baselineskip}
\item Since $G$ is connected, $\gamma_k>0$. For $\ell\ge k$,
Corollary~\ref{cor:cross-rank} gives
\[
 \frac{\delta_\ell}{\gamma_k}
 \ge\max\left\{\frac14,\frac{\ell-k+1}{2\ell}\right\}.
\]
For $\ell\ge2k-1$, the second term is strictly greater than $\frac{1}{4}$.

\item Suppose $\ell<k$. For $s\ge2$, let $H_s=K_{s,s}$. Then
\[
 A(H_s[K_2])=A(H_s)\otimes J_2+I_{2s}\otimes A(K_2),
\]
where $J_2$ is the all-one matrix.
The spectrum of $H_s[K_2]$ follows from
\cite[Corollary~2.2]{ACCV}; taking a Cartesian product adds the
eigenvalues of the second factor.

The Cartesian product
$G_s=H_s[K_2]\mathbin{\square}K_\ell$ is connected, nonbipartite and
vertex-transitive, of degree $2s+\ell$.
For large $s$, its spectrum gives
\[
 \delta_\ell=\frac{2\ell}{2s+\ell},\qquad
 \gamma_{\ell+1}=\frac{2s}{2s+\ell}.
\]
Since $k>\ell$,
$0\le \frac{\delta_\ell}{\gamma_k}\le\frac{\ell}{s}\to0$.

\item It remains to prove sharpness for $k\le\ell\le2k-2$.
Let $q\ge5$ be odd, and $r,s\ge1$ be integers. Take the vertex set
$\Z_q\times\{1,\ldots,s\}\times\{1,\ldots,r\}$.
For $1\le i\le r$, we call
$\Z_q\times\{1,\ldots,s\}\times\{i\}$ the $i$-th layer.
We define the adjacency matrix by
\[
A_{q,r,s}
 =A(C_q)\otimes J_s\otimes I_r
  +I_q\otimes I_s\otimes A(K_r),
\]
where $J_s$ is the all-one matrix. One observation is that the two terms in $A_{q,r,s}$ have disjoint support:
the first term changes the cycle coordinate, while the second term leaves that
coordinate fixed and changes the layer coordinate. Hence the matrix $A_{q,r,s}$ is
a simple graph adjacency matrix, and its degree is
\[
 d_{q,r,s}=2s+r-1.
\]
Within a fixed layer, consecutive cycle fibers are completely joined,
so that layer is connected. The second term joins the layers, and hence yields the
connectivity of the whole graph. Fixing a layer and a value of the second coordinate produces an odd cycle, so the graph is nonbipartite. The cycle rotations, permutations of the layers, and simultaneous permutations of the second
coordinate act transitively on its vertices.

The spectrum of $A(K_r)$ is $\{(r-1)^1,(-1)^{r-1}\}$ with multiplicities, the latter
part being absent when $r=1$. Taking tensor products of eigenvectors of
$A(C_q)$, $J_s$, and $A(K_r)$ gives an eigenbasis for $A_{q,r,s}$.
The eigenvalue $s$ of $J_s$ yields, for each integer $0\le t<q$,
the following normalized eigenvalues:
\begin{equation}\label{ed-eq:sharp-spectrum}
 \begin{aligned}
  &\frac{2s\cos\left(\frac{2\pi t}{q}\right)+r-1}{2s+r-1}
  &&\text{with multiplicity }1,\\
  &\frac{2s\cos\left(\frac{2\pi t}{q}\right)-1}{2s+r-1}
  &&\text{with multiplicity }r-1.
 \end{aligned}
\end{equation}
The eigenvalue $0$ of $J_s$, with multiplicity $s-1$, yields the remaining
normalized eigenvalues:
\[
\begin{aligned}
  &\frac{r-1}{2s+r-1}
  &&\text{with multiplicity }q(s-1),\\
  &\frac{-1}{2s+r-1}
  &&\text{with multiplicity }q(s-1)(r-1).
 \end{aligned}
\]

Note that terms with zero multiplicity are absent.

We first fix integers $q$ and $r$, and send $s\to\infty$. In
\eqref{ed-eq:sharp-spectrum}, exactly $r$ eigenvalues approach $1$ with $t=0$,
the next $2r$ eigenvalues approach $\cos\left(\frac{2\pi}{q}\right)$, and the lowest $2r$ ones approach
$-\cos\left(\frac{\pi}{q}\right)$. The additional eigenvalues from the zero eigenspace of $J_s$ tend uniformly to zero; the other cycle modes remain farther from the relevant spectral extremes.
Since $\cos\left(\frac{2\pi}{q}\right)>0$ and $\cos\left(\frac{\pi}{q}\right)>0$, these residual
eigenvalues lie strictly between the indicated upper and lower groups for
large $s$.

Fix $k\le\ell\le2k-2$ and set $r=k-1$. Then $r<k$ and $\ell\le2r$.
With $n=q rs$, as $s\to\infty$ we have
\[
 \mu_k\to\cos\left(\frac{2\pi}{q}\right),
 \qquad
 \mu_{n-\ell+1}\to-\cos\left(\frac{\pi}{q}\right).
\]
Consequently,
\[
 \lim_{s\to\infty}
 \frac{\delta_\ell}{\gamma_k}
 =\frac{1-\cos\left(\frac{\pi}{q}\right)}{1-\cos\left(\frac{2\pi}{q}\right)}
 =\frac{1}{4\cos^2\left(\frac{\pi}{2q}\right)}.
\]
Letting $q\to\infty$ through odd integers proves \eqref{ed-eq:sharp-gap}.\qedhere
\end{enumerate}
\end{proof}

\begin{example}[Sharpness of the partial-sum factor]\label{ex:prefix-sharp}
Fix $m\ge2$ and let $\ell\ge5$ be odd. Set
\[
 N_\ell=P(C_\ell)\oplus
 \begin{pmatrix}0&1\\1&0\end{pmatrix}^{\oplus(m-2)}.
\]
The second summand is absent when $m=2$. Here $\rho=1$.
The top $m$ eigenvalues consist of $m-1$ copies of $1$ and
$\cos\left(\frac{2\pi}{\ell}\right)$. The bottom $m$ eigenvalues consist of $m-2$
copies of $-1$ and two copies of $-\cos\left(\frac{\pi}{\ell}\right)$. Therefore
\[
 \sum_{j=1}^m\gamma_j=1-\cos\left(\frac{2\pi}{\ell}\right),\qquad
 \sum_{j=1}^m\delta_j=2\left[1-\cos\left(\frac{\pi}{\ell}\right)\right].
\]
Their ratio is
\[
 \frac{\sum_{j=1}^m\delta_j}{\sum_{j=1}^m\gamma_j}
 =\frac{1}{2\cos^2\left(\frac{\pi}{2\ell}\right)}\to\frac12.
\]
This proves the claimed sharpness at every fixed $m\ge2$.
\end{example}

Finally, we present an example to show the sharpness of the order of $d$ in the estimate given by Hu-Liu \cite[Theorem 1.3]{HL} for non-bipartite $d$-regular vertex-transitive graphs as follows:
\begin{equation}\label{eq:HL25}
\left[-1+c\frac{h_{out}^2}{d},1-C\frac{h_{out}^2}{d}\right]
\end{equation}
with some absolute constants $c$ and $C$. Note that all Cayley graphs are vertex-transitive, while the converse is not true.
\begin{example}[Sharpness of the order of $d$]
Fix an integer $m\ge 4$, let $G_m:=C(X_m,S_m)$
be the Cayley graph with respect to the group
$$X_m:=\Z/2\Z\times \Z/m\Z.$$
and the symmetric generating set
$$S_m:=\bigl(\{1\}\times \Z/m\Z\bigr)\cup\bigl(\{0\}\times\{1,-1\}\bigr).$$
Since $m\ge 4$, this set $S_m$ does not contain the identity.  Then $G_m$ is a simple $d_m$-regular abelian Cayley graph with degree
$$d_m=|S_m|=m+2.$$

Let
$$ L:=\{0\}\times\Z/m\Z, \qquad
  R:=\{1\}\times\Z/m\Z.$$
The generators in $\{1\}\times\Z/m\Z$ give all edges between $L$ and $R$.  Hence they form a complete bipartite graph $K_{m,m}$.  The two generators $(0,1)$ and $(0,-1)$ add a cycle $C_m$ inside $L$ and another cycle $C_m$ inside $R$.  Thus $G_m$ is obtained from $K_{m,m}$ by adding one cycle on each side. In other words, $G_m$ is the join of $m$-cycle $C_m$ with itself, i.e. $G_m=C_m \lor C_m$. In particular, $G_m$ is connected and non-bipartite. Then $G_m$ can achieve the lower bound in the estimate (\ref{eq:HL25}). We first need the following lemma for the eigenvalues of a join graph, also see in \cite[Exercise 1.13]{BH12}.
\begin{lemma}\label{lemma:join}
    Let $G_i=(V_i,E_i)~(i=1,2)$ be a $d_i$-regular~$(i=1,2)$ graph with $m_i=|V_i|~(i=1,2)$, and $H:=G_1\lor G_2$ be the join of $G_1$ and $G_2$. Assume that the eigenvalues of the (non-normalized) adjacency matrices of $G_1$ and $G_2$ are
    \[
\{\lambda_i\}_{i=1,2,\cdots,m_1-1}\cup \{\lambda_{m_1}=d_1\}
\]
    and
    \[
\{\mu_i\}_{i=1,2,\cdots,m_2-1}\cup\{\mu_{m_2}=d_2\},
\]
    respectively. Then the eigenvalues of the (non-normalized) adjacency matrix of $H$ are given by
    \[
 \text{Spec}(A(H))=
 \left\{
 \begin{aligned}
  &\lambda_1,\lambda_2,\cdots,\lambda_{m_1-1},
   \mu_1,\mu_2,\cdots,\mu_{m_2-1},\\
  &\frac{(d_1+d_2)\pm \sqrt{(d_1-d_2)^2+4m_1m_2}}{2}
 \end{aligned}
 \right\}.
\]
\end{lemma}
We know that it has two eigenvalues
\[
\lambda_-=-m+2,\qquad \lambda_+=m+2,
\]
and all other eigenvalues $\left\{2\cos \left( \frac{2k\pi}{m}\right)\right\}_{k=1,2,\cdots,m-1}$ lying in the interval $[-2,2]$. Hence the smallest eigenvalue of its normalized adjacency matrix is given by
\[
\mu_n(G_m)=\frac{-m+2}{m+2}=-1+\frac{4}{m+2}=-1+\frac{4}{d_m}.
\]
On the other hand, since $G
_m$ has a complete bipartite graph $K_{m,m}$ as its subgraph and the outer-vertex isoperimetry constant $h_{out}$ satisfies $h_{out}(K_{m,m})=1$, then we obtain that this constant on $G_m$ satisfies
\[
h_{out}(G_m) \ge h_{out}(K_{m,m})=1,
\]
while choosing one of the partition $L$ or $R$ yields the upper bound $h_{out}(G_m) \le1$. Hence we derive that
\[
h_{out}(G_m)=1.
\]
Therefore, it shows that the order $d^{-1}$ is optimal in the estimate (\ref{eq:HL25}). Particularly, an estimate with $d$ replaced by $d^\alpha$, where $\alpha<1$, cannot hold with a positive absolute constant. Indeed, since $h_{out}(G_m)=1$, such an estimate would force
\[
\frac{4}{d_m}\ge\frac{c}{d_m^\alpha}
\]
for all $m$, which is impossible as $m\to\infty$.
\end{example}

\appendix
\section{Even-walk-regularity of the five graph classes}
\label{app:five-classes}
In this appendix, we establish the even-walk-regularity of the five classes of connected nonbipartite graphs considered by Saha~\cite{Saha}, namely Cayley graphs, Cayley sum graphs, twisted Cayley graphs, twisted Cayley sum graphs, and vertex-transitive graphs.
\begin{lemma}\label{cg-lem:twosidedEWR}

Let $A$ be a real symmetric matrix indexed by a finite set $V$. Suppose a
group $H$ acts transitively on $V$ and, for each $h\in H$, there is a
permutation $\rho_h$ of $V$ such that
\[
 A_{x,y}=A_{hx,\rho_h y}\qquad(x,y\in V).
\]
Then $(A^{2t})_{xx}$ is independent of $x$ for every integer $t\ge0$.
\end{lemma}

\begin{proof}
Let $U_h$ and $V_h$ be the permutation matrices of $x\mapsto hx$ and
$y\mapsto\rho_h y$, respectively. The assumed identity gives $U_hA=AV_h$, and the symmetry of
$A$ also gives $V_hA=AU_h$. Hence
\[
 U_hA^2=AV_hA=A^2U_h.
\]
Every power $(A^2)^t$ therefore commutes with the transitive action of
$H$, thus its diagonal is constant.
\end{proof}

\begin{proposition}\label{cg-prop:fiveEWR}
Every finite simple graph in the following classes is even-walk-regular:
Cayley graphs, Cayley sum graphs, twisted Cayley graphs, twisted Cayley
sum graphs, and vertex-transitive graphs.
\end{proposition}

\begin{proof}
For a finite group $X$, a connection set $S\subseteq X$, and a bijection
$\sigma\colon X\to X$, consider the adjacency rules
\[
 y=\sigma(xs)\quad\text{and}\quad y=\sigma(x^{-1}s)
 \qquad(s\in S),
\]
The standard twisted Cayley and twisted Cayley sum graphs use these
rules with $\sigma$ an automorphism or anti-automorphism;
see~\cite{BSTwists}. The argument below permits
any bijection $\sigma$ whenever the resulting graph is simple.
The untwisted cases have $\sigma=\mathrm{id}$. Their
adjacency matrices satisfy
\[
 A^{\mathrm C}_{x,y}=\mathbf1_S(x^{-1}\sigma^{-1}(y)),
 \qquad
 A^{\Sigma}_{x,y}=\mathbf1_S(x\sigma^{-1}(y)).
\]
For each $g\in X$ and two vertices $x,y\in V$, the pairs of permutations $\tau_g$ and $\rho_g$ are given below,
\[
\begin{array}{c|cc}
 & \tau_g(x) & \rho_g(y)\\ \hline
 \text{twisted Cayley} & gx & \sigma(g\sigma^{-1}(y))\\
 \text{twisted Cayley sum} & xg^{-1} & \sigma(g\sigma^{-1}(y))
\end{array}
\]
and they satisfy $A_{x,y}=A_{\tau_g(x),\rho_g(y)}$. In each row the permutations
$\tau_g$ act transitively on $X$. Since the adjacency matrix is symmetric, Lemma~\ref{cg-lem:twosidedEWR} applies. Taking $\sigma=\mathrm{id}$ covers the ordinary Cayley and Cayley sum graphs. For a vertex-transitive graph, we take $H$ to be its automorphism group and set $\rho_h=h$.
\end{proof}

\section*{Acknowledgments}
This work was supported by the Scientific Research Innovation Capability Support Project for Young Faculty (Grant No. SRICSPYF-ZY2025160) and the National Natural Science Foundation of China (Grant No. 12431004). S. L. thanks Jack Koolen for raising the question of extending the results to walk-regular graphs following S. L.’s presentation of the work of Hu and Liu at the 2023 Workshop on Combinatorics and Graph Theory, hosted by Guangdong University of Technology.

\bibliographystyle{amsplain}
\renewcommand{\MR}[1]{}
\bibliography{references}
\end{document}